\documentclass[12pt]{article}
\usepackage{amssymb}
\title{On a variant of Hardy inequality between weighted Orlicz spaces}\author{Agnieszka Ka\l amajska \phantom{}\hskip 2mm
 and Katarzyna Pietruska-Pa\l{}uba\thanks{The work of both authors is
                         supported by a KBN grant no.
                         1-PO3A-008-29. It  is
also partially supported by EC FP6 Marie Curie programmes SPADE2
and CODY. } }

plus 6pt minus 12pt
\renewcommand{\it}{\sl}
\renewcommand{\em}{\sl}
\def\rn{{\mathbb{R}^{n}}}
\def\r{{\mathbb{R}}}
\newcommand{\vp}{\varphi}
\newcommand{\rp}{\mathbb{R}_+}
\newcommand{\ds}{\displaystyle}
\newcommand{\ts}{\textstyle}

\newtheorem{theo}{\bf Theorem}[section]
\newtheorem{coro}{\bf Corollary}[section]
\newtheorem{lem}{\bf Lemma}[section]
\newtheorem{rem}{\bf Remark}[section]

\newtheorem{ex}{\bf Example}[section]

\newtheorem{prop}{\bf Proposition}[section]

\makeatletter \@addtoreset{equation}{section}

\makeatother

\begin{document}
\maketitle \sloppy
\parindent 1em

\baselineskip=17pt

\begin{abstract}
{
 \protect\addtolength{\textwidth}{-4cm}
 \protect\addtolength{\oddsidemargin}{2cm}
Let $M$ be an $N$-function satisfying the $\Delta_2-$condition,
let $\omega, \vp$ be two other functions, $\omega\ge 0$. We study
Hardy-type inequalities
\[
 \int_{\rp} M(\omega (x)|u(x)|) {\rm exp}(-\vp (x))dx \le
C\int_{\rp} M(|u'(x)|) {\rm exp}(-\vp (x))dx ,
\]
where $u$ belongs to some  dilation invariant set ${\cal R }$
contained in the space of locally absolutely continuous functions.
We give sufficient conditions the triple $(\omega,\vp,M)$ must
satisfy in order to have such inequalities valid for $u$ from a
given set ${\cal R}$. The set  ${\cal R}$ can be smaller than the
set of Hardy transforms. Bounds for constants, retrieving
classical Hardy inequalities with best constants, are also given.
}
\end{abstract}



\section{Introduction}

{\bf General framework and classical approach.} Hardy type
inequalities have been  the subject of intensive research,  going
back to Hardy, who in the early 1920's (\cite{har20,hlp}) obtained
inequalities of the form
\begin{equation}\label{hhardy}
\int_{\rp} |u(t)|^pt^{\alpha -p}dt \le C\int_{\rp}
|u'(t)|^pt^{\alpha}dt.
\end{equation}

Over the decades to come, this inequality was generalized to
\[
\left( \int_{\rp} |\int_0^t f(\tau)d\tau |^qd\mu (t)
\right)^{\frac{1}{q}}\le C\left( \int_{\rp}|f(t)|^pd\nu (t)
\right)^{\frac{1}{p}},
\]
with two nonnegative Radon measures $\mu,\nu$, and further to
inequalities in the Orlicz-space setting:
\begin{equation}\label{ogoln}
Q^{-1}\left( \int_{\rp} Q(\theta (x)|Tf(x)|) w(x)dx \right)\le
P^{-1}\left( \int_{\rp} P(C\rho (x)|f(x)|) v (x)dx \right).
\end{equation}

Inequalities in  $L^p$ were studied  by Muckenhoupt \cite{muc},
Mazya and Rozin \cite{ma}, Bradley \cite{br}, Kokilashvili
\cite{ko}, Sinnamon \cite{sin}, Sawyer \cite{saw}, Bloom and
Kerman \cite{bloker}, Stepanow \cite{ste},  and many others (we
refer to the monographs \cite{kp,kufner,kuf3,kk,mipefi} and
references therein). As to the Orlicz-space result (\ref{ogoln}),
several authors contributed to the complete characterization of
admissible weights $\theta, w ,\rho, v$ and nondecreasing
functions $P,Q,$ which allow for (\ref{ogoln}) with
$Tf(x)=\int_0^x K(x,y) f(y)dy$ being the generalized Hardy
operator with kernel $K.$ To name just a few, we refer to the
papers of Bloom-Kerman \cite{bloker1,bloker}, Lai
\cite{lai1,lai2,lai3,lai4}, Heinig-Maligranda \cite{hm},
Bloom-Kerman \cite{bloker}, Heinig-Lai \cite{lai6} and  their
references.  For a more detailed account of such results, we refer
to  Sections \ref{lp1}, \ref{orliw1}.

Hardy-type inequalities are widely applicable in the PDE theory
and in functional analysis. For example, one can derive various
Sobolev embedding theorems in the $L^p$ setting, which can then be
 used to prove the  existence of solutions of the Cauchy
 problem in the elliptic and parabolic PDEs (see e.g.
 \cite{bct,gaap,kufner,cia,ma,st,kufap}), to study of the
 asymptotic behaviour of solutions
(\cite{aft,vz}), as well as their stability \cite{bv,cm}.  They
are present in the probability theory (see e.g. \cite{dy,rs,fi}).
Hardy inequalities are also of separate interest (see e.g.
\cite{kuf3,kp,tar}). For latest results, see the very recent
papers \cite{bct,chw,gop,ops} and their references.

The
 investigation of weighted or nonweighted
 Orlicz-Sobolev spaces defined by an $N-$function different from
 $\lambda^p$ is suggested by
physical models (see e.g. \cite{alber,ball,krzy2,po2,cia}).
Therefore it is worthwhile to examine Hardy-type inequalities in
general Orlicz spaces as well.

One of the  central problems evolving around Hardy inequalities
can be expressed as follows.
 Consider  sets of Hardy transforms ${\cal H}$ or conjugate Hardy transforms
 ${\cal H}^*$:

\begin{eqnarray}\nonumber
{\cal H}&=&\{ u(t)=\int_0^t f(s)ds,\ \int_0^a
|f(\tau)|d\tau <\infty\ \hbox{\rm for every}\ a>0  \} =\\
&=& \{ u\in W^{1,1}_{loc}(\rp) : \lim_{r\to 0} u(r)=0
\} ;\nonumber\\
 \nonumber {\cal H^*}&=&\{ u(t)=\int_t^\infty f(s)ds,\ \int_a^\infty |f(\tau)|d\tau <\infty\ \hbox{\rm
for every}\ a>0  \} =\\
&=& \{ u\in W^{1,1}_{loc}(\rp) : \lim_{r\to \infty} u(r)=0 \}
\label{hatra1} .
\end{eqnarray}

and let $Tf$ be either the Hardy transform of $f$,  or  the
conjugate Hardy transform of $f$. Then one deals with the
following problem.

{\sc Problem 1 (classical).} Given  two $N-$functions $Q,P$,
describe  all possible weights $(\theta,\omega,\rho,v)$ for which
inequality (\ref{ogoln}) follows  with some constant $C$
independent of $f$ (with $T$ the Hardy transform,  or the
conjugate Hardy transform).

This problem has been solved completely for the Hardy transform by
Bloom and Kerman in \cite{bloker} for modular functions $P$ and
$Q$ such that $Q$ dominates $P$ in some special sense (for the
details see Section \ref{orliw1}). Further generalization of
(\ref{lai}), without the domination restriction can be found in
Lai's paper \cite{lai4}. Therefore Problem 1 can be considered as
closed  and classical one.

{\bf Another approach and its motivation.} We are concerned with
another problem, which is expressed as follows.

{\sc Problem 2 (general).} Given  two $N-$functions $Q,P$  and
weights $(\theta,\omega,\rho,v)$, find a possibly big dilation
invariant set ${\cal R}$ contained is the set of locally
absolutely continuous functions
 for which the inequality
\begin{equation}\label{ogoln12}
Q^{-1}\left( \int_{\rp} Q(\theta (x)|u(x)|) w(x)dx \right)\le
P^{-1}\left( \int_{\rp} P(C\rho (x)|u'(x)|) v (x)dx \right) ,
\end{equation}
 would hold  with some constant $C$ independent of $u\in {\cal R}$.

Let us make several comments here.

It can happen that the given weights $(\theta,\omega,\rho,v)$ obey
the known requirements for the validity of (\ref{ogoln}), say for
$Tf$ being the Hardy transform of $f,$ as described in Problem 1.
In such a case the set ${\cal R}$ contains the full set of Hardy
transforms ${\cal H}$. On the other hand, if this requirement is
not fulfilled, we cannot expect (\ref{ogoln12}) to hold for every
$u\in {\cal H}$. In such a case  ${\cal R}\cap {\cal H}$ will be a
proper subset of ${\cal H}$. The solution to Problem 2 would
therefore lead to  Hardy-type inequalities within a (possibly)
narrower class of functions than this required so far.

Since the set of Hardy transforms (resp. conjugate Hardy
transforms) is dilation invariant, we want that  our set ${\cal
R}$ be dilation invariant as well. This  means by definition that
if $u\in {\cal R}$ then for every $\lambda
>0$ the function $u_\lambda (x):= u(\lambda x)$ also belongs to
${\cal R}$.



If we substitute $u=Tf$ in (\ref{ogoln12}), where $Tf$ is either
Hardy transform or conjugate Hardy transform, then $u'=f$, so that
(\ref{ogoln12}) reduces to the special case of (\ref{ogoln}).


{\bf The reduced problem and partial answers.} In this paper we
deal with the special variant of Problem 2, which reads as
follows.

{\sc Problem 3 (reduced).} Given  an $N-$function $M$ satisfying
the $\Delta_2-$condition, and a  pair of functions $(\omega ,\vp)$
where $\omega\ge 0$,  describe the possibly big dilation invariant
set ${\cal R}$ contained is the set of locally absolutely
continuous functions
 for which inequality
\begin{equation}\label{ogoln13}
 \int_{\rp} M(\omega (x)|u(x)|) {\rm exp}(-\vp (x))dx \le
C\int_{\rp} M(|u'(x)|) {\rm exp}(-\vp (x))dx ,
\end{equation}
 follows  with some constant $C$ independent on $u\in {\cal R}$.

 Variants of (\ref{ogoln13}) with $M(\lambda)=\lambda^p$ and
${\cal R}$ determined by the constraints concerning
$M,\omega,\phi$ have been studied e.g. in the papers
\cite{bee,co,flo,flowoj}, see also their references. To the best
of our knowledge their extension to Orlicz setting was not
considered so far.

Our main result formulated in Theorem \ref{drugie} states that if
the\ triple of functions $(\omega,\vp,M)$ satisfies certain simple
compatibility conditions, then we can indicate a dilation
invariant set ${\cal R}$ such that (\ref{ogoln13}) holds for every
$u\in {\cal R}$. Moreover, we  give some bounds on the constant
$C$, which can be expressed in terms of the  Simonenko lower and
upper index (see \cite{simon} and \cite{guspetre},
\cite{firokrbec} for interesting related results).

 The decision whether $u\in{\cal R}$ is based on its
behavior near zero and near infinity, which is very natural in
problems arising from PDE's: when analyzing a particular equation
one can often say that its solution (i.e. our function $u$) has
some `good' properties near the boundary, expressed in terms of
its rate of decay near the boundary.

As an illustration we derive the classical Hardy inequalities
(\ref{hhardy}) with best constants (see Section \ref{class01} for
discussion).

We also obtain sufficient conditions for (\ref{ogoln13}) to hold
for every $u\in {\cal H}$ (see Proposition \ref{osta-bloker} in
Section \ref{impli}). These conditions can be easily implemented
in practice, and since  the verification of the classical
Bloom-Kerman conditions (\ref{ewew}) from \cite{bloker} seems
rather hard, by using our approach one can avoid the verification
of Bloom-Kerman conditions and quickly deduce that (\ref{ogoln13})
is satisfied for every $u\in {\cal H}$.

Perhaps it is even more interesting to deal with the case when
Bloom-Kerman conditions (\ref{ewew}) are not satisfied, so that
inequality (\ref{ogoln13}) is not valid for every $u\in {\cal H
}$. Then we find a set ${\cal R}$ such that ${\cal R}\cap {\cal
H}$ is smaller than ${\cal H}$ and (\ref{ogoln13}) holds for every
its element (see Section \ref{impli1}).

As a particular type of inequalities alike (\ref{ogoln13}),  we
analyze those with $\omega = |\vp'|$ (see Section \ref{spec12}),
and  in the class of admissible $\vp'$s we obtain the inequality
\begin{equation}\label{caccio}
\int_{\rp}M(|\vp'(x)|u(x)|){\rm exp}(-\vp (x))dx \le C \int_{\rp}
M(|u'(x)|){\rm exp}(-\vp (x))dx.
\end{equation}
For $M(\lambda)=\lambda^p$  and $\psi(x)=
\exp\{-\frac{\vp(x)}{p}\}$ we get
\[
\int_{\rp}(|\psi'(x)u(x)|)^p dx \le C \int_{\rp} \left( |
u'(x)\psi(x)|\right)^p dx,
\]
which is nothing but a particular case of Caccioppoli inequality
on $\mathbb{R}^+$ (see e.g. \cite{cac,iwsbo}). Caccioppoli
inequalities are commonly used in the regularity theory, and so we
believe that our variant (\ref{caccio})  can be used in the
regularity theory as well.

\section{Preliminaries and  statements of  main
results}\label{prel}

\subsection{Preliminaries}

\noindent{\bf Orlicz spaces}

 Let us recall some preliminary facts
about Orlicz spaces, referring e.g. to \cite{rao-ren} for details.
Here we deal with Orlicz spaces of functions defined on $\rp$.

 Suppose that $\mu$ is a positive Radon
measure on $\rp$ and let $M:[0,\infty)\to[0,\infty)$ be an
$N-$function, i.e. a continuous convex function satisfying
$\lim_{\lambda\to 0}\frac{M(\lambda)}{\lambda}=0$ and
 \(\lim_{\lambda\to\infty}\frac{M(\lambda)}{\lambda}=\infty.\)

The weighted Orlicz space $L^{M}_\mu$ we deal with is by
definition the space
\[L^{M}_\mu\stackrel{def}{=}\{f:\rp\rightarrow\r\mbox{ measurable} :   \int_{\rp}
 M(\frac{|f(x)|}{K})d\mu(x)\leq 1\ \mbox{ for some }\ K >0\} ,  \]
 equipped with the Luxemburg norm
\[\|f\|_{L^M_\mu}=\inf\{ K>0 : \int_{\rp}M(\frac{|f(x)|}{K})d\mu(x)\leq
1\} .\] This norm is complete and turns $L^M_\mu$ into a Banach
space.  When $\mu$ is the Lebesgue measure, it is dropped from the
notation. For $M(\lambda)=\lambda^p$ with $p> 1$, the space
$L^M_\mu$ coincides with the usual $L^p_\mu$ space (defined on
$\rp$).

 The symbol $M^*$ denotes the complementary function of an
$N-$function  $M$, i.e its  Legendre transform: for $y,\geq 0,$
$M^*(y)=\sup_{x>0}[xy-M(x)].$  It is again an $N-$function and
from its definition we get the Young inequality:
\[
xy\leq M(x)+M^*(y), \;\;\mbox{ for } \;\; x,y\geq 0.
\]

 $M$ is said to
fulfill the $\Delta_2-$condition if and only if, for some constant
$c>0$ and every $\lambda >0,$ we have
\begin{equation}\label{delta2}
M(2\lambda)\leq cM(\lambda).
\end{equation}

In the class of differentiable convex functions the
$\Delta_2-$condition is equivalent to:
\[
\lambda M'(\lambda)\leq \tilde{c} M(\lambda),
\]
satisfied for every $\lambda >0$, with the constant $\tilde{c}$
being independent of $\lambda$ (see e.g. \cite{kr}, Theorem 4.1).

 We will need the following property of modulars:
 (see \cite{kr}, formula (9.21))
 \begin{equation}\label{norm1}
 \int_{\rp} M(\ts\frac{f(x)}{\|f\|_{L^M(\mu)}})\,d\mu(x)\leq
 1.
 \end{equation}
 When $M$ satisfies the $\Delta_2-$condition, then
 (\ref{norm1}) becomes an equality.

 The function $M_1$ is said to dominate
$M_2$ 
 if there exist two positive
constants $K_1,K_2$ s.t. $M_2(\lambda)\leq K_1 M_1(K_2\lambda)$
for every $\lambda>0.$  In such case we have
\begin{equation}
\|\cdot\|_{L^{M_2}_\mu}\leq K\|\cdot \|_{L^{M_1}_\mu},\;\; {\rm
with}\;\; K= K_2(K_1+1). \label{pstryk}\end{equation}
 Functions $M_1$
and $M_2$ are called equivalent 
when $M_2$ dominates $ M_1$ and $M_1$ dominates $ M_2$. In
particular equivalent N-functions give raise to equivalent
Luxemburg norms.

On the set ${\cal L}^M_\mu=\{u \mbox{ measurable:}\; \int
M(|u|)d\mu<\infty\},$ one introduces the so-called {\em dual}
norm:
\begin{equation}\label{dualnorm}
\|u\|_{L^{(M)}_\mu}={\rm sup} \{\int_{\rp} u(x)v(x)\,d\mu(x): v\in
L^{M^*}_\mu, \int_{\rp} M^*(|v(x)|)\,d\mu(x)\leq 1\}.
\end{equation}
The advantage of this norm is the H\"{o}lder-type inequality:
\begin{equation}\label{norms}
\mbox{for }\;\; f\in L^M_\mu, \;g\in L^{M^*}_\mu,\;\;\;\;
\int_{\rp} f\cdot g \,d\mu \leq \|f\|_{L^{(M)}_\mu}
\|g\|_{L^{(M^*)}_\mu},
\end{equation}

When $M$ satisfies the $\Delta_2-$condition, then ${\cal L}^M_\mu=
L^M_{\mu},$ and in general, the Orlicz space $L^M_\mu$ is the
completion of ${\cal L}^M_\mu$  in the dual norm.  The Luxemburg
norm and the dual norm are equivalent:

\[
\|u\|_{L^{(M)}_\mu}\leq \|u\|_{L^{M}_\mu} \leq
2\|u\|_{L^{(M)}_\mu}.
\]

{\bf\noindent Assumptions.} \\[1mm]
 Throughout the paper we assume:
\begin{description}
\item[(M)] $M:[0,\infty)\to [0,\infty)$ is a differentiable $N-$function, i.e.
 $M$ is convex, $M(0)=M'_+(0)=0$,
$M(\lambda)/\lambda \rightarrow \infty$ as $\lambda\to\infty$, and
moreover $M$ satisfies the condition:
\begin{equation}\label{pierwsze}d_M
\frac{M(\lambda)}{\lambda}\le M'(\lambda)\le D_M
\frac{M(\lambda)}{\lambda}\ \ \hbox{for every}\ \lambda >0,
\end{equation}
where $D_M\ge d_M\ge 1.$
\item [($\mu$)]$\mu$ is a  Radon measure on $\mathbb{R}_+,$ absolutely
continuous with respect to the Lebesgue measure and $\mu(dr)={\rm
e }^{-\vp(r)},$
 $\vp\in C^2(\mathbb{R}_+),$ $\vp'$ does never vanish,
 \item[($\omega$)] $\omega:(0,\infty)\to [0,\infty)$ is a $C^1-$function.
\end{description}

\begin{rem}\rm\label{zaba1}
The latter inequality in (\ref{pierwsze}) implies that $M$
satisfies $\Delta_2-$condition (see (\ref{delta2})). The condition
$d_M>1$ is equivalent to the fact that also $M^*$ satisfies the
$\Delta_2-$condition (see e.g. \cite{kr}, Theorem 4.3 or
\cite{akikppstud}, Proposition 4.1). Moreover, for any
$N-$function $M,$ the left-hand side in (\ref{pierwsze}) holds
with $d_M=1.$ \newline If $d_M$ and $D_M$ are the best possible
constants in (\ref{pierwsze}) they obey the definition of
Simonenko lower and upper index of $M$ and are related to Boyd
indices of $L^M(\rn,\mu)$ (see \cite{boy}, \cite{simon} for
definitions,  \cite{firokrbec}, \cite{guspetre}, \cite{yaq} for
discussion on those and other indices of Orlicz spaces).
\end{rem}

{\noindent\bf Conditions}\\[1mm]
 Before we proceed, we need to
introduce the following quantities.
We set
\begin{eqnarray}\label{ome}
\Omega &:=& \{ r\in \mathbb{R}_+:\omega(r)u(r)\neq 0\},\\
F&:= & \{ r\in \mathbb{R}_+:\omega(r)\neq 0,\
\omega'(r)\vp'(r)>0\},\nonumber\\
G&:= & \{ r\in \mathbb{R}_+:\omega(r)\neq 0,\
\omega'(r)\vp'(r)<0\},\nonumber
\end{eqnarray}
 Then we define:
\begin{eqnarray}
 b_1(r,\omega,\vp,M)&=& \left(
1+\frac{\vp''(r)}{(\vp'(r))^2} - \frac{\omega'(r)}{\omega
(r)\vp'(r)} \left[ d_M \chi_G(r) +D_M\chi_F(r) \right]
  \right);\nonumber
 \\[2mm]
b_1&=&b_1(\omega,\vp,M):= \inf \{ b_1(r,\omega,\vp,M) : r\in \mathbb{R}_+\}; \label{bjed} \\[2mm]
 b_2(r,\omega,\vp,M)&=& \left(
-1-\frac{\vp''(r)}{(\vp'(r))^2} + \frac{\omega'(r)}{\omega
(r)\vp'(r)}\left[ d_M\chi_F(r) + D_M\chi_G(r)\right]
  \right);\nonumber
 \\[2mm]
  b_2&=&b_2(\omega,\vp,M):= \inf \{ b_2(r,\omega,\vp,M) : r\in \mathbb{R}_+\}; \nonumber \\[2mm]
L=L(\omega,\vp) &:=&  \sup \left\{ \frac{\omega (r)}{|\vp'(r)|} :
r\in (0,\infty ) ,\vp'(r)\neq 0 \right\}\label{el0} .
\end{eqnarray}
We use the convention $\sup\emptyset =-\infty$, $\inf \emptyset
=+\infty$, $c/\infty =0$, $f\chi_A$ is the function $f$ extended
by $0$ outside $A$.

\subsection{Main results}

Our area of interest will be those triples $(M,\varphi,\omega)$
for which either {\small \begin{description}
\item[(B1)]
$b_1>0$, $L<\infty,
 $ or
\item[(B2)]
 $b_2>0$, $L<\infty$.
\end{description} }

We will deal with the  following function:
\begin{equation}\label{hr}
h^u(r)=h^{(u, \omega, \varphi, M)
}(r)=\frac{1}{\vp'(r)}M(\omega(r)|u(r)|),
\end{equation}
which is well defined since  $\varphi'(r)$ is never zero.

Let us introduce the following two classes of functions:
\begin{eqnarray}
{\cal R}^{+}_{(\omega,\vp, M)}&:=& \{ u \in W^{1,1}_{loc}(\rp):
\exists_{s_n\to 0, R_n\to\infty} :  \lim_{n\to \infty}\left(
h^u(R_n)\,{\rm e}^{-\vp (R_n)}-h^u(s_n){\rm e}^{-\vp (s_n)}\right)
\ge 0 \}; \nonumber \\
&&\label{admiplus}\\
{\cal R}^{-}_{(\omega,\vp, M)}&:=& \{ u \in W^{1,1}_{loc}(\rp):
\exists_{s_n\to 0, R_n\to\infty} :
  \lim_{n\to \infty}\left( h^u(R_n)\,{\rm e}^{-\vp
(R_n)}-h^u(s_n){\rm e}^{-\vp (s_n)}\right) \le 0 \}\nonumber;\\
&& \label{admiminus}
\end{eqnarray}
not precluding the limits  from being infinite. For simplicity we
will usually omit $(\omega,\varphi,M)$ from the notation.

Note that both sets ${\cal R}^+$ and ${\cal R}^-$ contain the set
of compactly supported $W^{1,1}$ functions and that they sum up to
the whole set $W^{1,1}_{loc}(\rp)$. Moreover, they are dilation
invariant, i. e. for every $\lambda >0$ and $u\in {\cal R}$ (where
${\cal R}$ is either ${\cal R}^+$ or ${\cal R}^-$) we have
$u_\lambda (x):= u(\lambda x)\in {\cal R}$.

\noindent Our main result reads as follows.

\begin{theo}\label{drugie}
Suppose that $M, \vp,\omega$ satisfying {\bf (M), ($\mu$),
($\omega$)} are such that Condition {\bf (B1)} (respectively {\bf
(B2)})  holds true. Then
\begin{equation}\label{hardy2}
  \int_{\rp} M (\omega (r))|u(r)|)\mu(dr)  \le C
\int_{\rp} M (|u'(r)|)\mu(dr)
\end{equation}
 holds for
every $u\in {\cal R}^{+}_{(\omega,\vp,M)}$ (resp. for every $u\in
{\cal R}^{-}_{(\omega,\vp,M)}$), where
$C=c(\frac{LD_M^2}{b_1d_M})$ (respectively
$C=c(\frac{LD_M^2}{b_2d_M})$), $c(x)={\rm max}(x^{d_M},x^{D_M})$.
 \end{theo}

As a direct consequence we also obtain the following theorem.

\begin{theo}\label{drugieab}
Suppose that $M, \vp,\omega$ satisfying {\bf (M), ($\mu$),
($\omega$)} are such that Condition {\bf (B1)} (respectively {\bf
(B2)})  holds true. Then the inequality
\[
  \|\omega u\|_{L^M_\mu}  \le \tilde{C}
\| u'\|_{L^M_\mu}
\]
 holds for
every $u\in {\cal R}^{+}_{(\omega,\vp,M)}$ (resp. for every $u\in
{\cal R}^{-}_{(\omega,\vp,M)}$),  where
$\tilde{C}=c(\frac{LD_M^2}{b_1d_M})+1$ (respectively
$\tilde{C}=c(\frac{LD_M^2}{b_2d_M})+1$), $c(x)={\rm
max}(x^{d_M},x^{D_M})$.
 \end{theo}

\section{Particular cases}

The main goal of this section is to illustrate Theorem
\ref{drugie} in various contexts.  At first we discuss the case
when $M(\lambda)=\lambda^p$ (Subsection \ref{lp1}). Then we turn
our attention to a general $M$ falling within our scope
(Subsection \ref{orliw1}).  Finally, in  Subsection \ref{spec12},
we restrict ourselves to the special choice of weights $\omega
=|\varphi'|$.

\subsection{Inequalities in the $L^p$ setting}\label{lp1}

When $M(\lambda)=\lambda^p$ with $p>1,$ our conditions get
simpler. In particular, $d_M=D_M=p,$ and since $\vp'$ is assumed
to be nonzero everywhere, we have:
\begin{eqnarray*}
b_1(r,\vp,\omega)&=&
1+\frac{\vp''(r)}{(\vp'(r))^2}-p\,\frac{\omega'(r)}{\omega(r)\vp'(r)}\chi_{\{\omega(r)\neq
0 \}}=-b_2(r,\vp,\omega).
\end{eqnarray*}
In particular, our Theorem yields results when $L<\infty,$ and
either $\inf_{r>0} b_1(r,\vp,\omega)>0,$ or $\sup_{r>0}
b_1(r,\vp,\omega)<0.$

\subsubsection{Classical Hardy inequalities}\label{class01}

As the first example illustrating our methods, we show that the we
can get the classical Hardy inequality.  This inequality
 reads as follows (see e.g. \cite{hlp}, Theorem 330 for
the classical source, Theorem 5.2 in \cite{kufner},  or \cite{kp}
for the statement,  historical framework and discussion).

\begin{theo}\label{hardyclas}
Let $1<p<\infty$, $\alpha\neq p-1$. Suppose that $u=u(t)$ is an
absolutely continuous function in $(0,\infty)$ such that
$\int_{0}^\infty |u'(t)|^p\,t^{\alpha}dt <\infty,$ and let
\begin{eqnarray*}
u^+(0)&:=& \lim_{t\to 0} u(t) =0\ \ {\rm for}\ \ \alpha <p-1,\\
u(\infty)&:=& \lim_{t\to\infty} u(t)=0\ \ {\rm for}\ \ \alpha
>p-1.
\end{eqnarray*}
Then the following inequality holds:
\begin{equation}\label{clasha}
\int_{0}^\infty |u(t)|^pt^{\alpha -p}dt \le C\int_0^\infty
|u'(t)|^pt^{\alpha}dt,
\end{equation}
where $C=\left( \frac{p}{|\alpha -p +1| }\right)^p$.
\end{theo}

\noindent We consider  the case $\alpha \neq 0$. Let us explain
how this theorem follows from our results.

Setting
\begin{eqnarray*}
M(r)=r^p, &   \mu (dr)=r^\alpha dr = {\rm exp}(\alpha \ln r), &
\omega(r)=\frac{1}{r}, \end{eqnarray*}
 we have \begin{eqnarray*} \vp
(r)=-\alpha\ln r, & \vp'(r)=\ds\frac{-\alpha}{r}, &
\vp''(r)=\frac{\alpha}{r^2},\\
  \omega'(r)=-\ds\frac{1}{r^2},&
\omega'(r)\vp'(r)=\frac{\alpha}{r^2}>0, & d_M=D_M=p.
\end{eqnarray*}
By a direct check we see that:
\begin{eqnarray*}
b_1= \frac{\alpha- (p-1)}{\alpha}, & b_2=
\ds\frac{(p-1)-\alpha}{\alpha}, &L=\frac{1}{|\alpha|}.
\end{eqnarray*}

Therefore for $\alpha >p-1$ and for $\alpha <0$ we have $b_1>0$,
while for $0<\alpha<p-1 $ we have $b_2>0$. In either case $\vp'$
does never vanish. The constant $C$ in (\ref{hardy2}) is equal to
$\left(\frac{p}{|\alpha-(p-1)|}\right)^p,$ which coincides with
the classical statement.

The only thing that remains to be checked is that any function $u$
as in the statement of Theorem \ref{hardyclas} for which the right
hand side in (\ref{clasha}) is finite belongs to ${\cal R}^+$ when
$\alpha
>p-1$ or $\alpha <0,$ and to ${\cal R}^-$ when $0 < \alpha
<(p-1).$  We will use standard arguments  (see e.g. \cite{kufner},
proof of Theorem 5.2).

 First suppose that $\alpha>p-1,$ and let $u$ be as in the assumptions of
 Theorem \ref{hardyclas}.  Then, for any
$t>0$ we define $U(t):=\int_t^\infty |u(\tau)|d\tau.$ One has:
\begin{eqnarray*}
U(t)&=&\int_t^\infty |u'(\tau)|\, d\tau \leq \int_t^\infty
|u'(\tau)|\tau^{\frac{\alpha}{p}}\tau^{-\frac{\alpha}{p}}d\tau
\\ &\le&
\left( \int_t^\infty |u'(\tau)|^p\tau^{{\alpha}}d\tau
\right)^{\frac{1}{p}} \left(\int_t^\infty
\tau^{-\frac{\alpha}{p-1}}d\tau\right)^{\frac{p-1}{p}}\\[2mm]
&=& \left( \frac{p-1}{\alpha -(p-1)} \right)^{\frac{p-1}{p}}\left(
\int_t^\infty |u'(\tau)|^p\tau^{{\alpha}}d\tau
\right)^{\frac{1}{p}} t^{-\frac{1}{p}(\alpha-(p-1))}<\infty.
\end{eqnarray*}

From this chain of inequalities and the condition $\int_0^\infty
|u'(\tau)|^p\tau^{{\alpha}}d\tau <\infty$ we infer not only that
$u\in W^{1,1}_{loc}(0,\infty)$ and that $U$ is well defined, but
also that $\lim_{R\to\infty}U(R)^pR^{\alpha-(p-1)}=0.$ Taking into
account that $\lim_{t\to\infty} u(t)=0,$ we have
\[
|u(t)|=|\int_t^\infty u'(\tau)d\tau |\le \int_t^\infty
|u'(\tau)|d\tau =U(t),
\]and therefore $\lim_{R\to\infty}|u(R)|^pR^{\alpha-(p-1)}=0$
as well.

Recall now the formulas defining the class ${\cal R}^+$.  The
function $h(r)e^{-\vp(r)}$ appearing there is in present situation
equal to $ -\frac{1}{\alpha}\,|u(r)|^pr^{\alpha-(p-1)},$ vanishing
for $r$'s tending to infinity. Therefore $u\in {\cal R}^+.$

When $\alpha <p-1,$ then we proceed similarly, but now we take
$U(t):=\int_0^t |u'(\tau)|d\tau.$ Again, $U$ is well defined and
$\lim_{r\to 0}U(r)^pr^{\alpha-(p-1)}=0,$ and since now $|u(r)|\leq
U(r)$ as well,  this permits to assert that $u\in {\cal R}^-.$ We
are done.

\subsubsection{General approach within $L^p$-spaces}

The following result is  considered classical now (see e.g.
\cite{ma}, Theorem 1 of Section 1.3.1).

\begin{theo}\label{mazja0}
Let $\mu,\nu$ be the nonnegative Borel measures on $(0,\infty)$,
let $\nu^*$ be the absolutely continuous part of $\nu$ and $1\le
p\le q\le \infty$. Then inequality
\begin{equation}\label{maz0}
\left( \int_{\rp} |\int_0^t f(\tau)d\tau |^qd\mu (t)
\right)^{\frac{1}{q}}\le C\left( \int_{\rp}|f(t)|^pd\nu (t)
\right)^{\frac{1}{p}}
\end{equation}
holds for an arbitrary  locally integrable function $f$ if and
only if
\[
B:= {\rm sup}\left(\mu [r,\infty)\right)^{\frac{1}{q}}\left(
\int_0^r (\frac{d\nu^*}{d\tau})^{-\frac{1}{p-1}}d\tau
\right)^{\frac{p-1}{p}} <\infty .
\]
\end{theo}

The case $p=q\ge 1$ is due to Muckenhoupt \cite{muc}.  Extensions
for general coefficients $p,q$ were proven by Mazya and Rozin
(\cite{ma}, Theorem 1 of Section 1.3.1), Bradley \cite{br} and
Kokilashvili \cite{ko}. Some other generalizations (admitting also
$p,q$ below 1) were obtained by Sinnamon \cite{sin}, Sawyer
\cite{saw}, Bloom and Kerman \cite{bloker}, Stepanow \cite{ste}
and others.

Observe that inequality (\ref{hardy2}) corresponding to
$\omega(r)=r$ is a particular case of (\ref{maz0}) when one takes
$\nu(dr)= e^{-\vp(r)}dr,$ $p=q,$ $\mu(dr)=r^pe^{-\vp(r)}dr$, but
only for the representant $u(t)=\int _0^t u'(\tau)d\tau$. In
particular $u^+(0)=0$, which we didn't require (Theorem
\ref{hardyclas}
 shows that in general the condition $u^+(0)=0$  {\em may not}
hold). In this case  (\ref{maz0}) reads
\[
 \int_{\rp} |u(\tau)\omega(\tau)|^pd\nu (\tau)  \le
C \int_{\rp}|u'(\tau)|^pd\nu (\tau) ,
\]
and $u(t)=\int_0^t u'(s)ds$, $u$ is an absolutely continuous
function.
 Condition $B<\infty$ is equivalent to
\[
 {\rm sup}_{r>0}\left(\int_r^{\infty}x^p{\rm exp }(-\vp
(x) )dx\right)\left( \int_0^r {\rm exp }(\frac{\vp (x)}{p-1} )dx
\right)^{p-1} <\infty .
\]
It is of different nature than our conditions {\bf (B1)} and {\bf
(B2)} and usually not  easy to handle. But since $B<\infty$ is
equivalent to the inequality (\ref{maz0}) holding for all $u$ in
the set of Hardy transforms ${\cal H}$ (see (\ref{hatra1})),
 our assumptions can serve as a tool towards
verifying $B<\infty.$

We may as well deal with the set of conjugate Hardy transforms
$u=-\int_t^\infty f(\tau)d\tau \in {\cal H}^*$ (see
(\ref{hatra1})) instead of $u=\int_0^t f(\tau)d\tau \in {\cal H}$
as illustrated in Theorem \ref{hardyclas}, the case $\alpha
>p-1$.

It may happen  that inequalities (\ref{maz0}) {\it do not hold} in
general on all set of Hardy transforms $u=\int_0^t f(\tau)d\tau\in
{\cal H}$, but they hold on a set essentially smaller than ${\cal
H}$, which we will indicate in the sequel. These are the sets
${\cal R}^+,$ ${\cal R}^-$ defined by (\ref{admiplus}) and
(\ref{admiminus}).

\subsection{Results in Orlicz spaces}\label{orliw1}

The series of papers of Mazya, Bloom-Kerman and Lai
\cite{ma,bloker,lai1,lai2,lai3,lai4,lai6} is concerned  with
inequalities in the form
\begin{equation}\label{lai}
Q^{-1}\left( \int_{\rp} Q(\omega (x)Tf(x)) r(x)dx \right)\le
P^{-1}\left( \int_{\rp} P(C\rho (x)f(x)) v (x)dx \right) ,
\end{equation}
where $Tf$ is the Hardy-type operator
\begin{eqnarray*}
Tf(x)=\int_0^x K(x,y) f(y)dy,  & x>0,
\end{eqnarray*}
with a suitable kernel $K.$ The case $K=1,$  corresponding to the
classical Hardy operator, is included there too. $P$ and $Q$ are
assumed to be certain nondecreasing functions on $\rp$ satisfying
\begin{eqnarray*}
\lim_{t\to 0^+} P(t)=\lim_{t\to 0^+} Q(t)=0, & \ds\lim_{t\to
\infty} P(t)=\lim_{t\to \infty} Q(t)=\infty .
\end{eqnarray*}
 Bloom and
Kerman proved in \cite{bloker} that within the class of modular
functions $P$ and $Q$ satisfying the following domination
property:
\begin{eqnarray*}
 \hbox{\it there exists a constant $\eta>0$ for which $ \sum
QP^{-1}(a_j)\le QP^{-1}(\eta\sum a_j), $}\\
\hbox{ whenever $\{ a_j\}$ is an arbitrary nonnegative
sequence},\end{eqnarray*} (\ref{lai}) is then equivalent to the
conditions
\begin{eqnarray}\label{blok}
\int_0^y P^* \left( \frac{G(\epsilon, y)K(y,x)}{B\epsilon v(x)\rho
(x) } \right)v(x)dx \le G(\epsilon ,y) <\infty\ \ {\rm and}\\
\int_0^y P^* \left( \frac{H(\epsilon, y)}{B\epsilon v(x)\rho (x) }
\right)v(x)dx \le H(\epsilon ,y) <\infty,\nonumber
\end{eqnarray}
holding for all $y>0$ and $\epsilon >0$, where $P^*$ is the
Legendre transform of $P$,
\begin{eqnarray*}
G(\epsilon, y)=PQ^{-1}\left( \int_y^\infty Q(\epsilon\omega
(x))r(x)dx\right),\\
H(\epsilon, y)=PQ^{-1}\left( \int_y^\infty Q(\epsilon\omega
(x)K(x,y))r(x)dx\right),
\end{eqnarray*}
$B>0$ is a constant.

In our particular case: $P=Q (=M)$ with $M$ satisfying the
$\Delta_2-$condition, $r=v={\rm exp}(-\vp (x))$, $K\equiv 1$,
$\rho\equiv 1$ inequality (\ref{lai}) reduces to
\begin{equation}\label{prot21}
 \int_{\rp} M(\omega (x)Tf(x)) {\rm exp}(-\vp (x))dx \le  C\left(
\int_{\rp} M(f(x)) {\rm exp}(-\vp (x)) dx \right),
\end{equation}
which is the type of inequality we are dealing with.

In this case conditions (\ref{blok}) are simpler and become
\begin{eqnarray}\label{blok1}
\int_0^y M^* \left( \frac{G(\epsilon, y)}{B\epsilon {\rm exp}(-\vp
(x))
 } \right){\rm exp}(-\vp (x))dx \le G(\epsilon ,y) <\infty ,
\end{eqnarray}
holding for all $y>0$ and $\epsilon >0, $ where
\begin{eqnarray*}
G(\epsilon, y)=\int_y^\infty M(\epsilon\omega (x)){\rm exp}(-\vp
(x))dx,
\end{eqnarray*}
$B>0$ is a constant. Further generalization of (\ref{lai}),
without the restriction $P<<Q,$ can be found in  Lai's  paper
\cite{lai4}.

Condition (\ref{blok1}) as well as Lai's condition are not
implemented easily: in practice, for given $\omega, \vp, M,$ it is
usually hard te see whether (\ref{blok1}) holds or not. Conditions
{\bf (B1)} and {\bf (B2)} are much simpler. When they are
satisfied, and when we know that the Hardy operator $Tf(x)$ (or
the dual Hardy operator $T^*f(x)=\int_x^\infty f(\tau)d\tau$, see
\cite{lai4}, the last remark on page 671) belongs to the set
${\cal R}^{-}$ or ${\cal R}^{+},$  then inequality (\ref{prot21})
is just the statement of Theorem \ref{drugie}.

\subsection{Special choice of weights. The case of $\omega
=|\vp'|$}\label{spec12}

Another case that substantially simplifies the approach is that of
$\omega=|\vp'|$ (in fact this was used in the proof of the
classical Hardy inequality). Since we require $\vp$ to be $C^1$
with nonzero derivative, $\vp'$ is either always positive, or
always negative.

This time around, we have:
\begin{eqnarray*}
b_1(r,|\vp'|,\vp,M)&=& \left\{
\begin{array}{ccc}
1+(1-d_M)\frac{\vp''(r)}{\vp'(r)^2} & {\rm if} &
\vp''(r)\le 0,\\[1mm]
1+(1-D_M)\frac{\vp''(r)}{\vp'(r)^2} & {\rm if} & \vp''(r)\ge 0,
\end{array}
\right.
\\[3mm]
b_2(r,|\vp'|,\vp,M)&=& \left\{
\begin{array}{ccc}
-1+(D_M-1)\frac{\vp''(r)}{(\vp'(r))^2} & {\rm if} &
\vp''(r)\le 0,\\[1mm]
-1+(d_M-1)\frac{\vp''(r)}{(\vp'(r))^2} & {\rm if} & \vp''(r)\ge 0,
\end{array}
\right.
\\[2mm]
L=1.&&
\end{eqnarray*}

\noindent In particular (as $D_M\ge d_M\ge 1$) we get:
\begin{eqnarray*}
b_1>0 && \mbox{ if and only if}\ \  \ {\rm sup}_{r>0}
\frac{\vp''(r)}{(\vp'(r))^2}< \frac{1}{D_M-1}=:\widetilde{D}_M;\\
 b_2>0 && \mbox{ if and only if}\ \ \
{\rm inf}_{r>0}\frac{\vp''(r)}{(\vp'(r))^2}>
\frac{1}{d_M-1}=:\widetilde{d}_M.
\end{eqnarray*}
This leads to the following conclusion, which is of separate
interest.

\begin{coro}\label{spac120}
Assume that conditions {\bf (M)}, {\bf ($\mu$)} are satisfied, and
\begin{eqnarray*}
{\rm either}\ \left( {\rm sup}_{r>0}\frac{\vp''(r)}{(\vp'(r))^2}<
\widetilde{D}_M \right) \ {\rm or}\ \ \left( {\rm inf}_{r>0}\;
\frac{\vp''(r)}{(\vp'(r))^2}> \widetilde{d}_M\right) .
\end{eqnarray*}
Then the inequalities
\begin{eqnarray*}
\int_{\rp} M(|\vp'(r)||u(r)|) {\rm exp}(-\vp (r))dr &\le&
C\int_{\rp} M (|u'(r)|){\rm exp}(-\vp (r))dr;\nonumber\\
\| \vp'u\|_{L^M_\mu} &\le& \tilde{C}\,  \|
u'\|_{L^M_\mu}
\end{eqnarray*}
hold for every $u\in {\cal R}^{+}_{(|\vp'| ,\vp, M )}$ with
$C=c(\frac{D_M^2}{b_1d_M})$, $\tilde{C}
=c(\frac{D_M^2}{b_1d_M})+1$ in the first case, and for every $u\in
{\cal R}^{-}_{(|\vp'| ,\vp, M )}$ with
$C=c(\frac{D_M^2}{b_2d_M})$, $\tilde{C}
=c(\frac{D_M^2}{b_2d_M})+1$ in the other case, $c(r)={\rm
max}(r^{d_M},r^{D_M})$.
\end{coro}

\begin{ex}[classical inequalities]\label{ilu1}\rm
To illustrate this corollary we consider again
\begin{eqnarray*}
M(\lambda)=\lambda^p, & p>1, & \varphi (\lambda)=-\alpha \ln r
\end{eqnarray*}
 as in Section \ref{class01}. In this case we have
 \begin{eqnarray*}
d_M=D_M=p, & \frac{\varphi''}{(\varphi')^2}\equiv
\frac{1}{\alpha}.
\end{eqnarray*}
We have
\begin{eqnarray*}
\frac{1}{\alpha}\equiv {\rm sup}_{r>0}\frac{\varphi'' (r)
}{(\varphi' (r))^2} < \frac{1}{p-1} =\tilde{D}_M & {\rm for} &
\alpha \in (-\infty, 0)\cup
(p-1,\infty),\\
\frac{1}{\alpha}\equiv {\rm inf}_{r>0}\frac{\varphi'' (r)
}{(\varphi' (r))^2} >  \frac{1}{p-1} =\tilde{d}_M & {\rm for}&
\alpha \in (0,p-1).
\end{eqnarray*}

 Therefore classical Hardy
inequality follows from Corollary \ref{spac120} as well.
\end{ex}

\begin{ex}[Hardy inequalities with logarithmic-type
weights]\label{logary}\rm As another illustration we show what
kind of inequality can be obtained for measures
$\mu(dr)=r^\alpha(\ln(1+r))^\beta \,dr,$ with $\alpha>0, \beta>0.$

In this case we have
\begin{eqnarray*}
\vp(r)&=:&\vp_{\alpha,\beta}(r)=-\alpha\ln(r)-\beta\ln\ln(1+r),\\[1mm]
\vp'(r)&=& \ds
-\frac{\alpha}{r}-\frac{\beta}{\ln(1+r)}\frac{1}{1+r},\\[1mm]
\vp''(r)&=& \ds\frac{\alpha}{r^2}+\frac{\beta}{(1+r)^2}
\frac{1}{\ln(1+r)}\left(1+\frac{1}{\ln(1+r)}\right).
\end{eqnarray*}
Choose $\omega(r)=|\vp'(r)|,$ then
\[b_1(r)=1+(1-D_M)\frac{\vp''(r)}{(\vp'(r))^2},\]
and so $b_1=1-(D_M-1)\sup_{r>0}\frac{\vp''(r)}{(\vp'(r))^2}.$ As
the derivative $\vp'(r)$ is of order $\frac{1}{r}$ and $\vp''(r)$
is of order $\frac{1}{r^2}$ on $\mathbb{R}_+,$ the supremum
involved is finite. Denote
\begin{equation}\label{sup}
s_{\alpha,\beta}=\sup_{r>0}\frac{\vp''(r)}{(\vp'(r))^2},
\end{equation}
then \begin{eqnarray*}
 b_1>0& \Leftrightarrow& \ds D_M < 1+\frac{1}{s_{\alpha,\beta}}.
 \end{eqnarray*}
As
$\frac{1}{r}\sim \omega(r),$  we arrive at the following.

\begin{theo}\label{logar}
Suppose $\alpha,\beta>0$ and let $M$ be such an $N-$function
satisfying {\bf (M)} that $D_M < 1+\frac{1}{s_{\alpha,\beta}},$
where $s_{\alpha,\beta}$ is given by (\ref{sup}). Then there
exists a constant $C>0$ such that the inequality
\[\int_{\mathbb{R}_+}M(\frac{u(r)}{r})r^\alpha(\ln(1+r))^\beta\,dr\leq
C\int_{\mathbb{R}_+}M(|u'(r)|)r^\alpha(\ln(1+r))^\beta\,dr \]
holds for all $u\in{\cal R}^+(|\vp'|,\vp,M)\supset
C_0^1(\mathbb{R}_+).$
\end{theo}

Similar analysis can be performed for negative $\alpha$ or
$\beta.$
\end{ex}

Practical information how to easily verify the assumption whether
$u\in {\cal R}^+$ or $u\in {\cal R}^-$ will be provided in
Subsection \ref{rylin}.

\section{Proofs of Theorems \ref{drugie} and \ref{drugieab}}

Before we pass to the actual proofs, let us formulate two easy
lemmas concerning Young functions. Although these lemmas may be
known to the specialists  (see e.g. \cite{guspetre}, \cite{simon}
and \cite{firokrbec} for related results), for readers'
convenience we submit the proofs.
\begin{lem}\label{lemm}
Suppose that $M$ is a differentiable $N-$function. Then we have.
\begin{description}
\item[(i)] Suppose that there exists a constant $D_M\geq 1$ such that
\begin{equation}\label{eqlemm1}
M'(r)\leq D_M\,\frac{M(r)}{r},\;\;\; \mbox{for every }\ r>0.
\end{equation}
Then for all $r>0, $ $\lambda \geq 1$
\[
M(\lambda r) \leq \lambda^{D_M}M(r).
\]
\item[(ii)]
Suppose that there exists a constant $d_M\geq 1$ such that
\[
 d_M\,\frac{M(r)}{r}\leq M'(r)\;\;\; \mbox{for every }\ r>0.
\]
Then for all $r>0, $ $\lambda \leq 1$
\[
M(\lambda r) \leq \lambda^{d_M}M(r).
\]
\item[iii)] Suppose that there exist constants $1\le d_M\le D_M$ such that
\[
 d_M\,\frac{M(r)}{r}\leq M'(r) \leq D_M\,\frac{M(r)}{r} \;\;\; \mbox{for every }\ r>0.
\]
Then for all $r>0, $ $\lambda >0$
\begin{equation}\label{ddwamodified}
M(\lambda r) \leq {\rm max}(\lambda^{d_M}, \lambda^{D_M}) M(r):=
c(\lambda) M(r).
\end{equation}
\end{description}
\end{lem}
 We recall Remark \ref{zaba1} about interpretation of constants
$d_M,D_M$.

\smallskip
 {\noindent\bf Proof.} We only prove part {\bf i)}.Part {\bf ii)}
 is proven analogously, while part {\bf iii)} is their direct consequence.\\
From (\ref{eqlemm1}) we get \( \frac{M'(r)}{M(r)}\leq
\frac{D_M}{r},\) and further, for any $r>0, \lambda>1$
\[\int_r^{\lambda r}\frac {M'(t)}{M(t)}\,dt\leq \int_r^{\lambda r}
\frac{D_M}{t}\,dt,\]
 which after integrating gives \([\ln M(t)]_r^{\lambda r}\leq[\ln
 t^{D_M}]_r^{\lambda r},\) and further \(M(\lambda r) \leq \lambda^{D_M}M(r).\)
 \hfill$\Box$

 \begin{lem}\label{lemmm}
 Suppose that $M$ is a differentiable $N-$function, and let $1\leq d_M\leq D_M$
 be two constants such that
 \begin{equation}\label{eqlemmm1}
d_M\,\frac{M(r)}{r} \leq M'(r)\leq D_M\,\frac{M(r)}{r},\;\;\;
\mbox{for every}\ r>0.
 \end{equation}
Then for every $r,s>0$ the following estimate holds true:
\begin{equation}\label{eqlemmm2}
\frac{M(r)}{r}\,s\leq \frac{D_M-1}{d_M}M(r)+\frac{1}{d_M}M(s).
\end{equation}
\end{lem}

 {\noindent\bf Proof.}
 Using the Young inequality $rs\leq M^*(r)+M(s)$ together with (\ref{eqlemmm1})
 we have:
 \begin{equation}\label{eqlemmm3}
\frac{M(r)}{r}\,s\leq\frac{1}{d_M}M'(r)s\leq
\frac{M^*(M'(r))}{d_M}+\frac{M(s)}{d_M}.
 \end{equation}
 From the very definition of the conjugate function $M^*$ we have $$M(r)=rM'(r)-M^*(M'(r)),$$
 and so
 \[
M^*(M'(r)) \leq D_M M(r) - M(r) =(D_M-1)M(r).
 \]
 Inserting this into (\ref{eqlemmm3}) we get
 (\ref{eqlemmm2}).\hfill$\Box$

{\bf Proof of Theorem \ref{drugie}.} Suppose   $u\in {\cal
R}^{+}_{(\omega,\vp,M)}$ (resp.
 $u\in {\cal R}^{-}_{(\omega,\vp,M)}$). Choose $s_n\to 0,
R_n\to\infty$ as in  (\ref{admiplus}) (resp. (\ref{admiminus})).
To abbreviate,  we denote
\begin{eqnarray*}
J=\int_{0}^\infty M (\omega (r))|u(r)|)e^{-\vp(r)}dr ,&
J_n=\ds\int_{s_n}^{R_n} M (\omega (r))|u(r)|)e^{-\vp(r)}dr ,\\
H=\ds\int_{0}^\infty M (|u'(r)|)e^{-\vp(r)}dr, &
H_n=\ds\int_{s_n}^{R_n} M (|u'(r)|)e^{-\vp(r)}dr,
\end{eqnarray*}

Let $h^u$ be given by (\ref{hr}). Under our assumptions, it is
well defined for every $r>0$. Since $u$ is $W^{1,1}_{loc}$ and $M$
is locally Lipschitz, we infer that $h^u\in W^{1,1}_{loc}
(\mathbb{R}_+)$,
\begin{eqnarray}\label{pochodna}
{(h^u)}'(r)&=&\frac{d}{dr}\left( \frac{1}{\vp'(r)} \right) M
(\omega (r)|u(r)|)\nonumber \\
&&+ \frac{1}{\vp'(r)}  M' (\omega (r)|u(r)|)\left(\omega'(r)|u(r)|
+\omega (r)u'(r){\rm sgn}\,u(r)\right),\nonumber\\
\phantom{}
\end{eqnarray}
in the sense of distributions and  almost everywhere,
 and $h$ is absolutely continuous on each interval
$[s,R]\subseteq (0,\infty)$
 (see e.g. \cite{ma}, Theorems 1 and 2, Sec. 1.1.3).
 Moreover,
 for every $R,s$ such that $0<s<R<\infty$

\begin{eqnarray*}
\int_s^R {(h^u)}'(r) e^{-\vp(r)}dr &=& h^u(r)e^{-\vp(r)}|_{s}^{R}
+\int_s^{R}M (\omega
(r)|u(r)|)e^{-{\vp(r)}} dr=\\
&:=&\theta (R,s)  +\int_s^{R}M (\omega (r)|u(r)|)e^{-{\vp(r)}} dr,
\end{eqnarray*}
and so we have
\begin{equation}\label{eee}
J_n= \int_{s_n}^{R_n} {(h^u)}'(r) e^{-\vp(r)}dr -\theta_n,
\end{equation}
where $\theta_n =\theta (R_n,s_n)\to\alpha\in [0, \infty] $ (resp.
$[-\infty,0]$ ).

 Inserting (\ref{pochodna}) inside (\ref{eee}) and using (\ref{ome})
yields, after some
 rearrangement,
\begin{eqnarray}
J_n&= & \int_{s_n}^{R_n}
(-\frac{\vp''(r)}{(\vp'(r))^2})M(\omega(r)|u(r)|)e^{-\vp(r)}dr \nonumber\\
&&-\int_{s_n}^{R_n} \frac{1}{|\vp'(r)|} M'(\omega(r)|u(r)|)
|\omega'(r)||u(r)|\chi_G(r)e^{-\vp(r)}dr \nonumber
\\
&&+\int_{s_n}^{R_n} \frac{1}{|\vp'(r)|}M'(\omega(r)|u(r)|)\left(
 |\omega'(r)||u(r)|\chi_F(r)+\omega
 (r)u'(r)\mbox{sgn}\, u(r)\right)e^{-\vp(r)}dr-\theta_n \nonumber\\
 &=& I_n-II_n +III_n -\theta_n.\label{truf}
 \end{eqnarray}

From now on, the proofs for the two cases: {\bf (B1)}, with $u\in
{\cal R}^+ $,
 and {\bf (B2)}, with $u\in
{\cal R}^-,$ differ slightly.\\

\textsc{Case 1.} Assumption {\bf (B1)} is satisfied,  $u\in
{\cal R}^+ $.\\

As $M$ satisfies  $-M'(\lambda)\le
-d_M\frac{M(\lambda)}{\lambda}$, one has
\begin{eqnarray}
-II_n &\le&-\int_{s_n}^{R_n} d_M M(\omega
(r)|u(r)|)\frac{|\omega'(r)|}{\omega(r)|\vp'(r)|}
\chi_{G}(r)e^{-\vp(r)}dr\nonumber\\
&=&\int_{s_n}^{R_n} d_M M(\omega
(r)|u(r)|)\frac{\omega'(r)}{\omega(r)\vp'(r)}
\chi_{G}(r)e^{-\vp(r)}dr =:IV_n\nonumber\\
\phantom{}\label{czteryn}
\end{eqnarray}
As to the estimate of $III_n,$ first estimate $ u'(r)\mbox{sgn}\,
u(r)$ by $|u'(r)|,$ and then use the inequality
 $M'(\lambda)\leq D_M\frac{M(\lambda)}{\lambda}.$  For every $r\in
 \Omega$ (see (\ref{ome})) we have:
\begin{eqnarray}
&&M'(\omega(r)|u(r)|)\left(|\omega'(r)|\,|u(r)|)\chi_{F }(r)+\omega(r)u'(r)\mbox{sgn}\,u'(r)\right) \nonumber \\
&\leq&
D_M\left[M(\omega(r)|u(r)|)\frac{|\omega'(r)|}{\omega(r)}\chi_{
F}(r) +
M(\omega(r)|u(r)|)\frac{|u'(r)|}{|u(r)|}\right]\nonumber\\[2mm]
 &=:&
D_M[A_1(r)+A_2(r)].\label{czapla}
\end{eqnarray}
This  implies
\begin{equation}\label{trzyn}
III_n\leq  {D_M} \int_{\Omega\cap [s_n,R_n]}
\frac{1}{|\vp'(r)|}[A_1(r)+A_2(r)]e^{-\vp(r)}dr=: V_n+VI_n.
\end{equation}
where
\begin{eqnarray*}
V_n&=& \int_{\Omega\cap [{s_n},{R_n}]} D_M
\frac{|\omega'(r)|}{\omega (r)|\vp'(r)|}\chi_{ F }(r) M(\omega
(r)|u(r)|)e^{-\vp(r)}dr\nonumber \\
&=& \int_{\Omega\cap [{s_n},{R_n}]} D_M \frac{\omega'(r)}{\omega
(r)\vp'(r)}\chi_{ F }(r) M(\omega
(r)|u(r)|)e^{-\vp(r)}dr,\nonumber \\
VI_n&=& \int_{\Omega\cap [s_n,R_n]}D_M
\frac{|u'(r)|}{|u(r)||\vp'(r)|} M(\omega (r)|u(r)|)e^{-\vp(r)}dr.
\end{eqnarray*}
Definition of the constant $b_1, $ (\ref{bjed}),
 yields
\[
b_1J_n\le J_n-I_n-IV_n-V_n,\]   and the series of estimates
(\ref{truf}), (\ref{czteryn}), (\ref{trzyn})  gives
\[
J_n= I_n-II_n+III_n-\theta_n\le I_n+IV_n+V_n+VI_n-\theta_n.\]
Combining the two we get
\[
b_1J_n\le VI_n-{\theta_n}.
\]
Consequently, since  $b_1$ is assumed to be positive,
\begin{equation}\label{jdm1}
J_n\leq\frac{D_M}{b_1}\int_{\Omega\cap [s_n,R_n]}
\frac{A_2(r)}{|\vp'(r)|}\, e^{-\vp(r)}dr-
\frac{\theta_n}{b_1}.\end{equation}

Now we use the estimates from Lemmas \ref{lemm} and \ref{lemmm}.
For any $0<\delta \leq 1$ (or, any $\delta >0$ when
$M(\lambda)=\lambda^p$) and $r\in\Omega$
\begin{eqnarray*}
A_2(r)&=& \delta
\omega(r)\frac{M(\omega(r)|u(r)|)}{\omega(r)|u(r)|}\,\frac{|u'(r)|}{\delta}\\[2mm]
&\leq&  \delta \omega(r)\left[
\frac{D_M-1}{d_M}\,M(\omega(r)|u(r)|)+\frac{1}{d_M}\,M(\frac{|u'(r)|}{\delta})\right]
\\[2mm]
&\leq &
  \delta \omega(r)\left[
\frac{D_M-1}{d_M}\,M(\omega(r)|u(r)|)+\frac{1}{d_M}c(\frac{1}{\delta})\,M(|u'(r)|)\right]
\end{eqnarray*}
(Lemma \ref{lemm} was used in the very last line).

 Using
now this estimate  on $A_2(r),$ we obtain from (\ref{jdm1})
\begin{eqnarray*}
J_n&\le&\frac{D_M\delta}{b_1}\left(\frac{D_M-1}{d_M}\,\int_{s_n}^{R_n}
\frac{\omega(r)}{|\vp'(r)|} M(\omega (r))|u(r)|)e^{-\vp(r)}dr
+\right.\\ &&+\left. \frac{1}{d_M}c(\frac{1}{\delta})
\int_{s_n}^{R_n} \frac{\omega (r)}{|\vp'(r)|}
M(|u'(r)|)e^{-\vp(r)}dr  \right)-\frac{\theta_n}{b_1}\\
&\le& \frac{D_M L\delta}{b_1}\left( \frac{D_M-1}{d_M}\, J_n+
\frac{1}{d_M }\,c(\frac{1}{\delta})H_n \right)
-\frac{\theta_n}{b_1},
\end{eqnarray*}
or, after rearranging :
 \[
\frac{\theta_n}{b_1}+\left(1-\frac{D_ML\delta}{b_1}\,
\frac{(D_M-1)}{d_M}\right)J_n\leq\frac{D_M L
\delta}{b_1d_M}\,c(\frac{1}{\delta})H_n.
\]
Choose now $\delta_0=\frac{b_1d_M}{LD_M^2},$  obtaining
 \begin{eqnarray*}
\frac{D_M\theta_n}{b_1}+J_n\leq c(\frac{1}{\delta_0})H_n.
\end{eqnarray*}

 Only now we  let $n\to\infty$. As
all limits: $\lim_{n\to\infty}J_n=J$, $\lim_{n\to\infty}H_n=H$,
$\lim_{n\to\infty}{\theta_n}$ are well defined and nonnegative
(finite or not), this implies
\[J\leq c(\frac{LD_M^2}{b_1d_M})H
\]
and finishes the proof.

 {\sc Case 2. } We now prove the statement under the assumption {\bf
 (B2)}, for $u\in {\cal R}^-$.
 We start with an expression similar to (\ref{truf}), but
 the integrals are  rearranged somewhat differently.

This time around we write
\begin{eqnarray}
J_n&= & \int_{s_n}^{R_n}
(-\frac{\vp''(r)}{(\vp'(r))^2})M(\omega(r)|u(r)|)e^{-\vp(r)}dr \nonumber\\
&&+\int_{s_n}^{R_n} \frac{1}{|\vp'(r)|} M'(\omega(r)|u(r)|)
|\omega'(r)||u(r)|\chi_Fe^{-\vp(r)}dr \nonumber
\\
&&-\int_{s_n}^{R_n} \frac{1}{|\vp'(r)|}M'(\omega(r)|u(r)|) \left(
|\omega'(r) ||u(r)|\chi_G +\omega
 (r)u'(r)\mbox{sgn}\, u(r) \right)e^{-\vp(r)}dr-\theta_n
\nonumber\\
 &=:& I_n+II_n' -III_n'-\theta_n.\label{truf1}
 \end{eqnarray}

Similarly as before, since $M'(\lambda)\ge
d_M\frac{M(\lambda)}{\lambda},$ we get
\begin{eqnarray*}
II_n' &\ge&\int_{s_n}^{R_n} d_M M(\omega
(r)|u(r)|)\frac{\omega'(r)}{\omega(r)\vp'(r)} \chi_{F
}(r)e^{-\vp(r)}dr=:IV_n'
\end{eqnarray*}
(one can omit the absolute values because $r\in F$). Further,
since $-M'(\lambda)\ge -D_M\frac{M(\lambda)}{\lambda}$,  for every
$r\in \Omega$ we have:
\begin{eqnarray*}
-M'(\omega(r)|u(r)|)\left(|\omega'(r)|\,|u(r)|)\chi_{G }(r)+\omega(r)|u'(r)|\right)  \\
\geq
-D_M\left[M(\omega(r)|u(r)|)\frac{|\omega'(r)|}{\omega(r)}\chi_{G
}(r)+
M(\omega(r)|u(r)|)\frac{|u'(r)|}{|u(r)|}\right]\\[2mm]
 :=
-D_M[B(r)+A_2(r)],
\end{eqnarray*}
where $A_2(r)$ is the same as in (\ref{czapla}).
 This last estimate combined with (\ref{truf1})
imply
\[
-III_n'\geq -{D_M} \int_{\Omega\cap [s_n,R_n]}
\frac{1}{|\vp'(r)|}[B(r)+A_2(r)]e^{-\vp(r)}dr:= V_n'+VI_n'.
\]
On the other hand,
\[
b_2J_n\le I_n +IV_n' +V_n'-J_n\ {\rm and}\ J_n\ge
I_n+II_n'-III_n'-\theta_n\ge I_n +IV_n' +V_n' +VI_n'-\theta_n.
\]
We get $ b_2J_n\le -VI_n' +\theta_n$, which leads to
\[
J_n\leq\frac{D_M}{b_2}\int_{\Omega\cap [s_n,R_n]}
\frac{A_2(r)}{|\vp'(r)|}\, e^{-\vp(r)}dr+ \frac{\theta_n}{b_2}.\]
Now the proof follows along the same lines as the proof in the
first case, starting from (\ref{jdm1}) up to its end, with $b_2$
replacing $b_1$ and $\theta_n$ replacing $-\theta_n$.  We are
done. \hfill$\Box$

\begin{rem}\rm
 Observe that  for
$M(\lambda)=\lambda^p$ we have $d_M=D_M=p,$ $c(x)=x^p.$ Therefore
under the assumptions of Theorem \ref{drugie} the constant $C$
equals to either $\left(\frac{Lp}{b_1}\right)^p$ if $b_1>0$ or $
\left(\frac{Lp}{b_2}\right)^p$ if $b_2>0$.
\end{rem}

{\bf Proof of Theorem \ref{drugieab}:} Without loss of generality
we can assume that $A=\| u'\|_{L^M_\mu}$ is finite. Let us
substitute $u_A:= \frac{u}{A}$ in (\ref{hardy2}). Then we get
\[
\int_{\rp} \frac{M}{C}\left(\frac{\omega |u|}{A}\right){\rm
exp}(-\vp)\le 1,
\]
which implies $\| \omega u\|_{L^{\frac{M}{C}}_\mu}\le A$. As
$N-$functions $\frac{M}{C}$ and $M$ are equivalent, we get by
(\ref{pstryk}) that $\| \omega\mu\|_{L^M_\mu} \le (C+1)\| \omega
u\|_{L^{\frac{M}{C}}_\mu}\le (C+1)A =(C+1)\| u'\|_{L^M_\mu}$.
Therefore the result follows.\hfill$\Box$

\begin{rem}\rm One may compare our results with those recently proven in
\cite{bct}. Namely, in Theorem 3.1 p. 416  the authors obtain the
following inequality for the Gaussian measure, for
$M(\lambda)=\lambda^p:$
\[
\int_{\rp} |\omega u|^p {\rm exp}(-\frac{x^2}{2})dx \le
(\frac{p}{p-1})^p \int_{\rp}|u'|^p{\rm exp}(-\frac{x^2}{2})dx,
\]
where $\omega (x)=\frac{{\rm exp}(x^2/2(p-1))}{\int_0^x {\rm
exp}(\sigma^2/2(p-1))d\sigma}$, holding for every $u\in
W_0^{1,p}(\rp , d\mu)$ (i.e. the completion of $C_0^\infty (\rp) $
in the weighted Sobolev  space $W^{1,p}(\rp , d\mu)$, where $\mu
(x)={\rm exp}(-x^2/2)dx$ is the Gaussian measure). In this  case
we have $\vp (x)=x^2/2$ and the quantity $\omega(x)/\vp'(x)$ is
not bounded as  required by our Theorem 3.1. Instead, the weight
$\omega$ obeys a different requirement, which is the ODE: $x\omega
- (p-1)\omega'=(p-1)\omega^2$.
\end{rem}

\section{Analysis of sets ${\cal R}$}\label{ry}

\subsection{Verification of the condition $u\in {\cal R}$}\label{rylin}

The analysis provided in this subsection is two-fold. At first we
show an easy practical method to verify  whether $u\in {\cal R}$
(here ${\cal R}$ equals to ${\cal R}^+$ or ${\cal R}^-$). Further
analysis is devoted to the discussion when the condition
\begin{eqnarray}
\int_{\rp} M(|u'(r)|){\rm exp}(-\vp (r))dr <\infty\label{0drugie}.
\end{eqnarray}
together with $u\in {\cal H}$ (respectively $u\in {\cal H}^*$)
implies that $u\in {\cal R}$.

We start with the following result.

\begin{prop}\label{sprawdzanie}
Suppose $M,\vp$ satisfy conditions ${\bf (M)}$, {\bf ($\mu$)
}{\bf($\omega$)}, $\vp\in C^2(\rp)$, $\vp'$ does never vanish and
$L<\infty$ (see \ref{el0}). Then the following statements hold
true.
\begin{description}
\item[(i)]
Assume that $\vp'(R)\to 0$ as $R\to\infty,$ and $u\in
W^{1,1}_{loc}(\rp)$ be such that $u(R)$ together with
$u(R)e^{-\vp(R)}$ are bounded next to $\infty$. Then for $\vp'<0$
we have $u\in {\cal R}^{+}_{(\omega,\vp, M)}$, while for $\vp'>0$
we have $u\in {\cal R}^{-}_{(\omega,\vp, M)}$
\item[(ii)]
Assume that $\vp'(r)\to 0$ as $r\to 0,$ and  and $u\in
W^{1,1}_{loc}(\rp)$ be such that $u(r)$ together with
$u(r)e^{-\vp(r)}$ are bounded next to $0$. Then for $\vp'>0$ we
have $u\in {\cal R}^{+}_{(\omega,\vp , M)}$, while for $\vp'<0$ we
have $u\in {\cal R}^{-}_{(\omega,\vp, M)}$.
\end{description}
\end{prop}

{\noindent \bf Proof.} {\bf i)}: We have ($h^u$ was defined by
(\ref{hr})):
\begin{eqnarray*} |h^u
(R)|{\rm exp}(-\vp (R)) &=&\frac{M (\omega
(R)|u(R)|)}{|\vp'(R)}{\rm exp}(-\vp (R)) \\
&=&\frac{M (|\vp'(R)||u(R)|\frac{\omega
(R)}{|\vp'(R)|})}{|\vp'(R)|}{\rm exp}(-\vp (R)) \\&\le &c(L)
\frac{M (|\vp'(R)||u(R)|)}{|\vp'(R)|}{\rm exp}(-\vp (R))\chi_{\{
u(R)\neq 0\} },
\end{eqnarray*}
where $c(\cdot)$ is defined in (\ref{ddwamodified}). Property
$\frac{M(r)}{r}\to 0$ as $r\to 0$ and our assumptions imply
\[
\frac{M(|\vp'(R)||u(R)|)}{|\vp'(R)||u(R)|}\chi_{\{ u(R)\neq 0\}
}\to 0\ \ \hbox{\rm as}\ \  R\to\infty\] and so $c(L)|u(R)|{\rm
exp}(-\vp (R))$ is bounded next to $\infty$. Therefore $h^u
(R)|{\rm exp}(-\vp (R))\to 0$ as $R\to\infty$ and the statement
now follows from definition of ${\cal R}^{+}$ and ${\cal R}^{-}$.

{\bf ii)}: Similarly as in the proof of part {\bf i)}, we check
that $|h^u (r)|{\rm exp}(-\vp (r))\to 0$ as $r\to 0$.
 \hfill$\Box$

\vspace{2em}
 In the remaining part of this subsection we examine the
property (\ref{0drugie}).

\vspace{1em}

To proceed, let us set some additional notation.\\

First, recall $c(\cdot)$ from (\ref{ddwamodified}), and then
define for $r>0$
\begin{eqnarray}\label{defgie}
 f_{\vp}(r)= c^{-1}(e^{-\vp(r)}),\\
 A_{\vp}(r)=\|\frac{1}{f_\vp}\|_{L^{(M^*)}(0,r)}\;\;\;\;
B_{\vp}(r)=\|\frac{1}{f_\vp}\|_{L^{(M^*)}(r,\infty)}.\nonumber
\end{eqnarray}
The norms considered here are the dual norms defined by
(\ref{dualnorm}).

We will distinguish two naturally appearing cases: the first  one
when $A_\vp$ is well defined, and the other when $B_\vp$ is well
defined.\\

{\bf\noindent CASE 1.} $A_\vp(r)$ well defined for small $r'$s.

We are now to analyze the case when $u\in {\cal H}$ which
satisfies  (\ref{0drugie}) belongs to ${\cal R}$.

  We start with a
following lemma.

\begin{lem}\label{lemlem}
Assume that\\
1. $M,\vp, \omega$ satisfy {{\bf (M), ($\vp$), ($\omega$)}},\\
2. $A_\vp$ is well defined for small $r$'s
and  the function
\[
K(r)=\frac{ M(\omega (r) A_\vp(r))}{|\vp'(r)|}\,e^{-\vp (r)}
\]
is bounded next to 0.

Then for every  $u\in W^{1,1}_{loc}(\rp)$  such  that
$\int_0^\infty M(|u'(r)|)e^{-\vp(r)}\,dr<\infty$\  and \\
$\lim_{r\to 0} u(r)=0$ the function $h^u$ defined by (\ref{hr})
satisfies
\begin{equation}\label{eqeq}
\lim_{r\to 0} h^u(r)e^{-\vp(r)}=0.
\end{equation}
\end{lem}

{\bf\noindent Proof.} Set
\[U(r)=\int_0^r |u'(\rho)|d\rho.\]
From inequality (\ref{norms}) we have
\begin{eqnarray}
 \int_0^r|u'(\rho)|d\rho & =& \ds\int_0^r |u'(\rho) f_\vp (\rho)|\cdot
\frac{1}{|f_\vp (\rho)|}d\rho\leq \|u'\cdot f_\vp
\|_{L^{(M)}(0,r)}\cdot A_{\vp}(r).\nonumber\\
\phantom{}\label{abojawiem}
 \end{eqnarray}
We have
\begin{equation}\label{kwiatek}
\| u'f_\vp\|_{L^M (0,r)} \le \| u'\|_{L^M_\mu (0,r)}.
\end{equation}

 Indeed, if $A= \| u'f_\vp\|_{L^M(0,r)}$ then we get from
 the definition of $f_\vp$ (\ref{defgie}) and the property (\ref{ddwamodified}):
\[
1=\int_{(0,r)} M(\frac{|u'|f_{\vp}}{A})dx \le \int_{(0,r)} c(f_\vp
) M(\frac{|u'|}{A})dx =\int_{(0,r)} M(\frac{|u'|}{A}) {\rm
exp}(-\vp)dx.
\]
Therefore $A\le \| u'\|_{L^M_\mu (0,r)}$, which gives
(\ref{kwiatek}).

As  dual and Luxemburg norms are equivalent, we get
\[ \|u' f\|_{L^{(M)}(0,r)}\leq A\|u'\|_{L^{(M)}_{\mu}
(0,r)}\leq A
 \|u'\|_{L^{(M)}_{\mu} (0,\infty)}<\infty ,\]
 where $A$ is some universal constant.
Therefore $U$ is well-defined.

From the assumption $\lim_{r\to 0} u(r)=0$ we get $|u(r)|\leq
U(r),$ and the estimate (\ref{abojawiem}) holds true for $u$
instead of $U$ as well.

 Now,
\begin{eqnarray*}
|h^u(r)e^{-\vp(r)}|&=&\ds
\frac{M(\omega(r)|u(r)|)}{|\vp'(r)|}e^{-\vp(r)}\leq
\frac{M(\omega(r)\|u'\|_{L^{(M)}_\mu(0,r)}A_{\vp}(r))}{|\vp'(r)|}\,e^{-\vp(r)}\\[2mm]
&\leq & \frac{M(\omega(r)
A_{\vp}(r))}{|\vp'(r)|}\,e^{-\vp(r)}c(\|u'\|_{L^{(M)}_\mu(0,r)})=K(r)c(\|u'\|_{L^{(M)}_\mu(0,r)}).
\end{eqnarray*}

But since $c(x)\to 0$ when $x\to 0,$ and
$\|u'\|_{L^{(M)}_\mu(0,r)}\to 0$
 when $r\to 0,$ the assertion (\ref{eqeq}) follows from the
 boundedness of $K(r)$ for small $r'$s. \hfill$\Box$

As a corollary, we obtain straight from the definition of ${\cal
R}^+, {\cal R}^-$:
\begin{coro}\label{comcom}
Suppose that the assumptions 1, 2 of Lemma \ref{lemlem}
are satisfied. Then we have.
\begin{description}
\item[i)] When $\vp'>0,$ then  \[
\{u\in {\cal H}, \ \int_{\rp} M(|u'(r)|){\rm exp}(-\vp (r))dr
<\infty \} \subset {\cal R}^+ ;\]
\item[ii)] When $\vp'<0,$ then \[
\{u\in {\cal H},\ \int_{\rp} M(|u'(r)|){\rm exp}(-\vp (r))dr
<\infty\} \subset {\cal R}^-.\]
\end{description}
\end{coro}

To illustrate the  statements  above we now discuss the following
example.

\begin{ex}\rm
Let us consider the case $\varphi (r)=-\alpha\ln r$, $\alpha
<p-1$, $M(\lambda)=\lambda^p$ and $\omega (r)=\frac{1}{r}$ as in
Theorem \ref{hardyclas}. In such a case we have
$M^*(\lambda)=c_p\,\lambda^{\frac{p}{p-1}}$ and, in this range of
$\alpha'$s,
\begin{eqnarray*}
f_{\vp}(r)=r^{\frac{\alpha}{p}}, & M^{*}(f_{\vp}(r))=c_p
r^{\frac{\alpha}{p-1}}, & A_\vp(r)= a_p r^{\frac{p-1-\alpha}{p}}.
\end{eqnarray*}
Then $K(r)$ is just a constant. Moreover, every  function $u\in
{\cal H}$ which satisfies $\int_0^\infty
|u'|^px^{\alpha}dx<\infty$ belongs to ${\cal
R}^{-}_{(\frac{1}{x},x^{\alpha},\lambda^p)}$ for $\alpha <0,$ and
to ${\cal R}^+_{(\frac{1}{x},x^{\alpha},\lambda^p)}$ when
$\alpha>0.$
\end{ex}

{\bf \noindent CASE 2.} $B_\vp(R)$ well defined for large $R'$s.

In this case we have the following dual statements dealing with
the property that $u\in {\cal H}^*$ which satisfies
(\ref{0drugie}) belongs to ${\cal R}$.

\begin{lem}\label{lemlemlem}
Assume that\\
1. $M,\vp, \omega$ satisfy {{\bf (M), ($\vp$), ($\omega$)}},\\
2. $B_\vp$ is well defined
and  the function
\[
L(R)=\frac{ M(\omega (R) B_\vp(R))}{|\vp'(R)|}\,e^{-\vp (R)},
\]
is bounded next to infinity.

Then for every $u\in {\cal H}^*$ is such  that $\int_0^\infty
M(|u'(r)|)e^{-\vp(r)}\,dr<\infty$ the function $h^u$ defined by
(\ref{hr}) satisfies
\[
\lim_{R\to \infty} h^u(R)e^{-\vp(R)}=0.
\]
\end{lem}

{\noindent \bf Proof.} It is almost identical with that of Lemma
\ref{lemlem}. We must replace now $U(r)=\int_0^r |u'(\rho)|d\rho$
with $U^*(R)=\int_R^\infty |u'(\rho)|d\rho$ and then proceed as
before.\hfill$\Box$

As a counterpart of Corollary \ref{comcom} we assert the
following.

\begin{coro}\label{comcomcom}
Suppose that the assumptions of Lemma \ref{lemlemlem}
are satisfied. Then we have.
\begin{description}
\item[i)] When $\vp'>0,$ then  \[
\{ u\in {\cal H}^*, \ \int_{\rp} M(|u'(r)|){\rm exp}(-\vp (r))dr
<\infty \} \subset {\cal R}^- ;\]
\item[ii)] When $\vp'<0,$ then \[
\{u\in {\cal H},\ \int_{\rp} M(|u'(r)|){\rm exp}(-\vp (r))dr
<\infty\} \subset {\cal R}^+.\]
\end{description}
\end{coro}

This result is illustrated by following example.

\begin{ex}\rm
We again consider the case $\varphi (r)=-\alpha\ln r$,
$M(\lambda)=\lambda^p$ and $\omega (r)=\frac{1}{r}$ as in Theorem
\ref{hardyclas}, but now $\alpha >p-1$. In such a case we have
$B_\vp (r)=b_p r^{\frac{-\alpha +p-1}{p}}$ and $L(r)$ is just a
constant. Therefore every  function $u\in {\cal H}^*$ which
satisfies $\int_0^\infty |u'|^px^{\alpha}dx<\infty$ belongs to
${\cal R}^{+}_{(\frac{1}{x},x^{\alpha},\lambda^p)}$.
\end{ex}

\begin{rem}\rm
Suppose that the assumptions of Theorem \ref{drugie} are satisfied
and let us  put ${\cal R}={\cal R}^+$ in case of {\bf (B1)} and
${\cal R}={\cal R}^-$ in case of {\bf (B2)}. One could ask whether
 the spaces
\[
\{ u\in {\cal R} :\int_{\rp} M(|u'(r)|){\rm exp}(-\vp (r))dr
<\infty \} ,
\]
playing the crucial role in the inequality (\ref{hardy2}) can
possibly be nonlinear. We do not know the answer to this question.
\end{rem}

\subsection{Relation to Bloom-Kerman results}\label{rbk}

We are now to compare our results with that of Bloom and Kerman
from \cite{bloker}.

\subsubsection{When our conditions imply Bloom and Kerman
ones}\label{impli}

Here we will  show an example where our assumptions yield
inequality (\ref{prot21}) for all functions $Tf= \int_0^t
f(\tau)d\tau \in{\cal H},$  so inherently the Bloom-Kerman
condition is satisfied.

\noindent Corollary \ref{comcom} results in the following
proposition.

\begin{prop}\label{osta-bloker}
Assume that $(M,\vp,\omega)$ satisfy the assumptions 1,2 in Lemma
\ref{lemlem} and either  $\left( \vp'>0, b_1>0, L<\infty \right)$
or $\left( \vp'<0, b_2>0, L<\infty \right)$. Then we have.
\begin{description}
\item[(i)] There exists a constant $C>0$ such that inequality
\[
  \int_{0}^\infty M (\omega (x)|u(x)|){\rm exp}(-\varphi
(x))dx  \le C \int_{0}^\infty M (|u'(x)|){\rm exp}(-\varphi (x))dx
\]
 holds for every    $u\in {\cal H}$.
\item[(ii)]
The triple $(\omega,\varphi ,M)$ satisfies the Bloom-Kerman
condition:
\begin{equation}\label{ewew} \int_0^y M^* \left( \frac{G(\epsilon,
y)}{B\epsilon {\rm exp}(-\varphi (x))
 } \right){\rm exp}(-\varphi
(x))dx \le G(\epsilon ,y) <\infty ,
\end{equation}
holding for all $y>0$ and $\epsilon >0$,
\[
G(\epsilon, y)=\int_y^\infty M(\epsilon\omega (x)){\rm
exp}(-\varphi (x))dx,
\]
$B>0$ is a constant.
\end{description}
\end{prop}

{\bf\noindent Proof.} {\bf i):}  Statement {\bf i)} is  just a
combination of Corollary \ref{comcom} and Theorem \ref{drugie}.

{\bf ii):} Let $u(x)=\int_0^x f(\tau)d\tau =(Tf)(x)\in {\cal H }$
be the Hardy transform of $f$. Then part {\bf i)} implies
\[
  \int_{0}^\infty M (\omega (x)Tf(x)){\rm exp}(-\varphi
(x))dx  \le C \int_{0}^\infty M (f(x)){\rm exp}(-\varphi (x))dx,
\]
 which is proven
to be equivalent to the condition (\ref{ewew}) (see Theorem 1.7 in
\cite{bloker} and our comments in Subsection
\ref{orliw1}).\hfill$\Box$

\subsubsection{When Bloom and Kerman conditions are not
satisfied}\label{impli1}

It may happen that our conditions are satisfied and the
Bloom-Kerman conditions are not. In such a case inequalities
(\ref{hardy2}) cannot hold for every Hardy transform $u\in{\cal
H}$ (see (\ref{hatra1})) but they hold on proper subsets in the
set of Hardy transforms.
This is illustrated by the following result.

\begin{prop}
There exists a triple $(\omega,\varphi,M)$ such that conditions
${\bf (M)}$,  ${\bf (\mu)}$, ${\bf (\omega)}$  are satisfied and
\begin{description}
\item[i)]
$(\omega,\varphi,M)$ satisfies {\bf (B1)}, in particular
inequality
\begin{equation}\label{tatatata}
  \int_{\rp} M (\omega (r))|u(r)|)\mu(dr)  \le C
\int_{\rp} M (|u'(r)|)\mu(dr)
\end{equation}
is satisfied whenever $u\in {\cal R}^{+}_{(\omega,\varphi, M)}$.
\item[ii)] $(\omega,\varphi,M)$ does not satisfy the Bloom-Kerman condition
(\ref{ewew}).
\item[iii)]
The set \[ {\cal R}^{(0,+)}_{(\omega,\varphi, M)}:= {\cal
R}^{+}_{(\omega,\varphi, M)}\cap {\cal H}\subseteq {\cal H}\] is a
proper subset of ${\cal H}$. Moreover, inequality (\ref{tatatata})
is not satisfied for every $u\in {\cal H}$, with the constant
independent of $u$.
\end{description}
\end{prop}

{\bf Proof.} Let $p>1$ and
\begin{eqnarray*}
M(\lambda)=\lambda^p, & \varphi(x)=-\frac{1}{2}x^2, & \omega
(x)=x.
\end{eqnarray*}
In particular conditions ${\bf (M)}$,  ${\bf (\mu)}$, ${\bf
(\omega)}$ are satisfied,
 $b_1 =1+\frac{1}{x^2}(p-1)>0$,
$L=1$ and the condition {\bf (B1)} is also satisfied. Therefore
the result {\bf i)} follows by Theorem \ref{drugie}.

\noindent {\bf ii):} The Bloom-Kerman condition does not hold: one
has just
\[G(\epsilon,y)= \int_y^\infty (\epsilon
x)^pe^{\frac{1}{2}x^2}dx=\infty,
\]
and so (\ref{ewew}) is violated.

\noindent {\bf iii):} The Laplace function
\[
 u(r)=\int_0^r {\rm exp}(-\tau^2)d\tau,
\]
belongs to ${\cal H}\setminus {\cal R}^{0,+}_{(\omega,\varphi,
M)}$ and does not satisfy (\ref{tatatata}).
\hfill$\Box$

 {\bf Acknowledgement.}
The authors would like to thank Miroslav Krbec  for helpful advice
and discussions.

{\small
 \vspace{1em}
 {\sc Address:} Institute of Mathematics,
 University of  Warsaw,
 ul.\ Banacha 2,
 02-097 Warszawa, Polska (Poland),
 emails: {\tt kalamajs@mimuw.edu.pl} and  {\tt kpp@mimuw.edu.pl}
}

\end{document}